\documentclass[12]{amsart}
\usepackage{amsfonts}
\usepackage{amssymb}
\newtheorem{defi}{Definition}[section]
\newtheorem{satz}[defi]{Theorem}
\newtheorem{lemm}[defi]{Lemma}

\newtheorem{cor}[defi]{Corollary}
\newtheorem{prop}[defi]{Proposition}

\newcommand{\supp}{\operatorname{supp}}

\newcommand{\trace}{\operatorname{trace}}

\newcommand{\h}{\mathcal{H}}
\newcommand{\HH}{\mathcal{H}}
\newcommand{\TT}{\mathbb{T}}

\newcommand{\lspace}{\vspace{0.2cm}}

\newcommand{\DD}{\mathcal{D}}

\newcommand{\C}{\mathbb{C}}

\newcommand{\T}{\mathbb{T}}
\newcommand{\N}{\mathbb{N}}
\newcommand{\R}{\mathbb{R}}

\newcommand{\LL}{\mathcal{L}}
\newcommand{\bmodd}{\operatorname{{\mathrm{BMO^d}}}}
\newcommand{\bmos}{\operatorname{{\mathrm{BMO_{so}}}}}

\newcommand{\bmosd}{\operatorname{{\mathrm{BMO_{so}^d}}}}

\newcommand{\bmow}{\operatorname{{\mathrm{WBMO^d(\LL(\h))}}}}
\newcommand{\bmoh}{\operatorname{{\mathrm{BMO^d(\LL(\h))}}}}

\newcommand{\bmond}{\operatorname{{\mathrm{BMO_{norm}^d}}}}

\newcommand{\bmop}{\operatorname{{\mathrm{BMO_{para}}}}}
\newcommand{\bmol}{\operatorname{{\mathrm{BMO_{mult}}}}}
\newcommand{\opf}{\operatorname{{\mathcal{F}_{00}}}}

\newcommand{\mat}{\mathrm{Mat}(\C, n \times n)}

\newcommand{\bmo}{\operatorname{{\mathrm{BMO}}}}
\newcommand{\sbmod}{\operatorname{{\mathrm{SBMO^d}}}}
\newcommand{\wbmod}{\operatorname{{\mathrm{WBMO^d}}}}

\newcommand{\eps}{\varepsilon}

\begin{document}

\title[Operator BMO spaces]
{Operator-valued dyadic BMO spaces}
\author{Oscar Blasco}
\address{Department of Mathematics,
Universitat de Valencia, Burjassot 46100 (Valencia)
 Spain}
\email{oscar.blasco@uv.es}
\author{Sandra Pott}
\address{Department of Mathematics, University of Glasgow,
University Gardens, Glasgow G12 8QW} \email{sp@maths.gla.ac.uk}
\keywords{operator BMO, Carleson measures}
\thanks{{\it 2000 Mathematical Subjects Classifications.}
                               Primary 42B30, 42B35, Secondary 47B35 \\
The first author gratefully acknowledges support by the LMS and
Proyectos MTM 2005-08350 and PR2006-0086. The second author
gratefully acknowledges support by EPSRC and by the Nuffield
Foundation.}
%\footnote{{\it 2000 Mathematical Subjects Classifications.}
 %                              Primary 42B30, 42B35, Secondary 47B35}
\begin{abstract}
We consider BMO spaces of operator-valued functions, among them
the space of operator-valued functions $B$
which define a bounded paraproduct on $L^2(\h)$. We obtain several
equivalent formulations of $\|\pi_B\|$ in terms of the norm of the
"sweep" function of $B$  or of averages of the norms of
martingales transforms of $B$ in related spaces. Furthermore, we
investigate a connection between John-Nirenberg type inequalities
and Carleson-type inequalities via a product formula for
paraproducts and deduce sharp dimensional estimates for John-Nirenberg
type inequalities.
\end{abstract}
\maketitle
\section{Introduction}
Spaces of BMO functions on the real numbers $\R$ or the circle $\T$,
taking values in the bounded linear operators on a Hilbert space,
 have been investigated in a number of different
contexts in recent years, for example non-commutative $L^p$ spaces
\cite{pxu}, \cite{mei1}, matrix-weighted inequalities
\cite{gptv2}, \cite{gptv},
 sharp estimates for vector Carleson Embedding Theorem
\cite{katz}, \cite{ntv}, \cite{nptv}, \cite{petermichl}, observation operators in linear systems over
 contractive semigroups
\cite{jp}, \cite{jpp1}, and Hankel operators in several variables
\cite{ps}.

The theory of operator valued BMO functions is much more complicated than the scalar theory
and remains to be fully understood.
Some of the different yet
equivalent characterizations of scalar $\bmo(\T)$ or $\bmo(\R)$ lead to distinct spaces of operator valued
BMO functions. In many cases, we can express this in the
language of operator spaces by saying that different operator space structures on the scalar
BMO space arise naturally from the different yet equivalent characterisations of scalar BMO.
These difficulties reflect partly the subtle geometric properties of the dual Banach
space $\LL(\HH)$ of bounded
linear operators on a Hilbert space.

It is often easier to consider dyadic versions of BMO and to work with dyadic versions of classical operators
like the Hilbert transform $H$ or the Hankel operator with symbol $b$, $\Gamma_b$. Two such dyadic
counterparts of a Hankel operator $\Gamma_b$
 are the \emph{dyadic paraproduct $\pi_b$} and the operator $\Lambda_b = \pi_b + \pi_{\bar b}^*$.
While the former has a natural interpretation as a Carleson Embedding operator, the latter connects
more easily in the operator valued case to the theory of vector-valued BMO functions
(in particular to the space $\bmond(\LL(\HH))$). Estimates for Hankel operators can then be
obtained by averaging techniques.

One important difference between the scalar-valued and the operator-valued settings is the failure
of a certain version of the classical \emph{John-Nirenberg Lemma}, or in other words, the lack of
boundedness
of the ``\emph{sweep}'', which governs the behaviour of the dyadic paraproduct.

The purpose of the present paper is  to study in particular the
spaces arising from the operators $\pi_b$ and $\Lambda_b$,
 to investigate the relationship between dyadic paraproduct, its ``real part'' $\Lambda_b$ and the sweep,
and to give sharp dimensional estimates for the sweep in the ``strong'' BMO norm $\|\cdot\|_{\bmosd}$ and other norms,
answering a question in \cite{gptv2}.

 Let
$\DD$ denote the collection of dyadic subintervals of the unit
circle $\T$, and let $(h_I)_{I \in \DD}$, where $h_I =
\frac{1}{|I|^{1/2}} ( \chi_{I^+} - \chi_{I^-})$, be the Haar basis
of $L^2(\T)$. Let $\h$ be a separable, finite or
infinite-dimensional Hilbert space and let  $\opf$ denote the
subspace of $\LL(\h)$-valued functions on $\T$ with finite formal
Haar expansion.
  Given $e,f\in \h $ and $B \in L^2(\T,\LL(\h))$, we denote by $B_e$ the
function in $L^2(\T,\h)$ defined by $B_e(t)= B(t)(e)$
and by
$B_{e,f}$ the function in $L^2(\T)$ defined by $B_{e,f}(t)=
\langle B(t)(e),f\rangle$.  As in the scalar case,
  let $B_I$ denote the formal Haar coefficients
$\int_I B(t) h_I dt$, and $m_I B = \frac{1}{|I|} \int_I B(t) dt$ denote the
average of $B$ over $I$ for any $I \in \DD$. Observe that for $B_I$
and $m_IB$ to be well-defined operators, we shall be assuming that
the $\LL(\h)$- valued function $B$ is $weak^*$-integrable. That
means, using the duality $\LL(\h)=(\h\hat\otimes \h)^*$, that
$\langle B(\cdot)(e),f\rangle\in L^1(\T)$ for $e,f\in \h $. In
particular, for any measurable set $A$, there exist $B_A\in
\LL(\h)$ such that
$\langle B_A(e),f\rangle=\langle\int_A B(t)(e) dt, f\rangle
$.

Let us denote by $\bmodd(\T,\HH)$ the space of Bochner integrable
$\h$-valued  functions $b: \T \rightarrow \h$ such that
\begin{equation}
   \|b\|_{\bmodd(\h)}=\sup_{I \in \DD} (\frac{1}{|I|} \int_I \| b(t) - m_I b\|^2
dt)^{1/2}<\infty
\end{equation}
and by $\wbmod(\T,\h)$ the space of Pettis integrable $\h$-valued
functions $b: \T \rightarrow \h$ such that
\begin{equation}
   \|b\|_{\wbmod(\h)}=\sup_{I \in \DD, e \in \HH, \|e\|=1} (\frac{1}{|I|} \int_I |\langle b(t) - m_I b,
e \rangle|^2 dt)^{1/2}<\infty
\end{equation}

 Let us define different version of dyadic operator-valued
BMO to be considered throughout the paper.

 We denote by $\bmond(\T,\LL(\h))$
the space of Bochner integrable $\LL(\h)$-valued  functions $B$
such that
\begin{equation} \label{bmond}
   \|B\|_{\bmond}=\sup_{I \in \DD} (\frac{1}{|I|} \int_I \| B(t) - m_I B
\|^2 dt)^{1/2}<\infty.
\end{equation}
 and  denote by ${\rm WBMO^d}(\T,\LL(\h))$
the space of $weak^*$-integrable $\LL(\h)$-valued  functions $B$
such that
\begin{multline}\label{wbmo}
\|B\|_{{\rm WBMO^d}}=\sup_{I \in \DD, \|e \|=\|f \|=1}
                    (\frac{1}{|I|} \int_I | \langle(B(t) - m_I B)e,f\rangle
|^2 dt)^{1/2} \\= \sup_{e \in \HH, \|e \|=1}
\|B_e\|_{\wbmod(\T,\HH)} < \infty,
\end{multline}
or, equivalently, such that
$$
\|B\|_{{\rm WBMO^d}}= \sup_{A \in S_1, \| A \|_1 \le 1 } \|
\langle B, A \rangle \|_{\bmodd(\T)} < \infty.
$$
Here, $S_1$ denotes the ideal of trace class operators in
$\LL(\h)$, and $\langle B, A \rangle$ stands for the scalar-valued
function given by $\langle B, A \rangle (t) = \trace( B(t) A^*)$.

In the operator-valued setting one has another natural
formulation.  Denote by ${\rm SBMO^d}(\T,\LL(\h))$ the space of
$\LL(\h)$-valued functions $B$ such that $B(\cdot)e\in
\bmodd(\T,\HH)$ for all $e\in\h$ and such that
\begin{equation}\label{sbmo}
  \|B\|_{{\rm SBMO^d}}= \sup_{I \in \DD,e \in \HH, \|e \|=1}
                    (\frac{1}{|I|} \int_I \| (B(t) - m_I B)e \|^2
                    dt)^{1/2}< \infty.
\end{equation}

We would like to point out that while $B$ belongs to one of the
spaces $ \bmond(\T,\LL(\h))$ or ${\rm WBMO^d}(\T,\LL(\h))$)  if
and only if $B^*$ does, this is  not the case for the space
$\mathrm{SBMO}^d(\T,\LL(\h))$. This leads to the following notion:

\begin{defi} (see \cite{gptv2}, \cite{petermichl},\cite{pxu} )  We say that
 $B\in\bmosd(\T,\LL(\h))$, \label{bmosd} if
 $B$ and $B^*$ belong to $ {\rm SBMO^d}(\T,\LL(\h))$.
We define $\|B\|_{\bmosd}= \|B\|_{{\rm SBMO^d}}+\|B^*\|_{{\rm
SBMO^d}}.$
\end{defi}

Continuous versions of this space in the more general setting of
functions taking values in a von Neumann algebra with a semifinite
normal faithful trace
 were studied by Pisier and Xu \cite{pxu} and more recently
by Mei \cite{mei1}, together with an $H^p$ theory and a rich
duality and interpolation theory.

\lspace

We now define another operator-valued BMO space, using the notion
of Haar multipliers. As in the scalar-valued case (see
\cite{per}), a sequence $(\Phi_I)_{I \in \DD}$, $\Phi_I\in
L^2(I,\LL(\h))$ for all $I\in \DD$, is said to be an
\emph{operator-valued Haar multiplier}, if there exists $C>0$ such
that
$$\|\sum_{I\in \DD}\Phi_I(f_I)h_I\|_{L^2(\T,\h)}\le C (\sum_{I\in
\DD}\|f_I\|^2)^{1/2} \text{ for all } (f_I)_{I \in \DD} \in
l^2(\DD,\h).$$ We write  $\|(\Phi_I)\|_{mult}$ for the norm of the
corresponding operator on $L^2(\T,\h)$.

Let us observe that \begin{equation} \label{estimacion}
\|\Phi_J\|_{L^2(\T,\h)}\le\|(\Phi_I)\|_{mult} |J|^{1/2}, \quad
J\in \DD.
\end{equation}
\begin{defi} Let us define $P_I B
=\sum_{J\subseteq I} h_JB_J,$  and  use the notation
$$\Lambda_B(f)=\sum_{I \in \DD} (P_I B) (f_I) h_I.$$ We define $ \bmol (\T,\LL(\h))$ as the
space of those $weak^*$-integrable $\LL(\h)$-valued functions for
which $(P_IB)_{I\in\DD}$ defines a bounded operator-valued Haar
multiplier,
    and write
 \begin{equation}\label{bmol}\|B\|_{\bmol}= \|\Lambda_B \|=
\|(P_IB)_{I\in\DD}\|_{mult}.
\end{equation}
\end{defi}

Let us now give the definition of a further BMO space, the space
defined in terms of dyadic paraproducts.

Let $B \in \opf$. We define
$$
   \pi_B: L^2(\T, \h) \rightarrow L^2(\T, \h), \quad f = \sum_{I \in \DD}
 f_I h_I\mapsto
                                 \sum_{I \in \DD}  B_I (m_I f) h_I,
$$
and
$$
   \Delta_B: L^2(\T, \h) \rightarrow L^2(\T, \h), \quad f = \sum_{I \in \DD}
 f_I h_I\mapsto
                                 \sum_{I \in \DD}  B_I (f_I)
\frac{\chi_I}{|I|}.
$$
$\pi_B$ is called the vector paraproduct with symbol $B$.

 It is  elementary to see that
\begin{equation}  \label{mult}
  \Lambda_B( f )=
                                 \sum_{I \in \DD}  B_I (m_I f) h_I
                            + \sum_{I \in \DD} B_I (f_I) \frac{\chi_I}{|I|}.
\end{equation}

This shows that $\Lambda_B=\pi_B+\Delta_B$. Observe  that
$\Delta_B=\pi^*_{B^*}$. Therefore $(\Lambda_B)^*=\Lambda_{B^*}$,
and $\|B\|_{\bmol}=\|B^*\|_{\bmol}$.

\begin{defi}
Let $E_k B = \sum_{|I| > 2^{-k} } B_I h_I$ for $k \in \N$. The space $
\bmop (\T,\LL(\h))$ consists  of those $weak^*$-integrable
operator-valued functions for which $\sup_{k \in \N} \|\pi_{E_k
B}\| < \infty$.
 For such functions,
$\pi_B f = \lim_{k \to \infty} \pi_{E_k B} f$ defines
a bounded linear operator on $L^2(\T,\h)$, and we write \begin{equation}
\label{bmop}\|B\|_{\bmop}=
\|\pi_B\|.
\end{equation}
\end{defi}

Let us notice that
\begin{equation}\label{form}
\Lambda_B f= B  f -\sum_{I\in \DD}(m_IB)(f_I) h_I.
\end{equation}
From here one concludes immediately that
\begin{equation}  \label{prop:linfty}
 L^\infty(\T,\LL(\h))\subseteq \bmol(\T,\LL(\h)).
\end{equation}
However, Tao Mei \cite{mei} has shown recently that $L^\infty(\TT, \LL(\HH)) \nsubseteq \bmop$
and therefore in particular $\bmol \nsubseteq \bmop$. This is in
contrast to the situation of scalar paraproducts in two variables,
where $\bmol(\TT^2) = \bmop(\TT^2)$ (\cite{BP}, Thm 2.8).

The following   chain of strict inclusions for
infinite-dimensional $\h$ can be shown (see \cite{new}):
\begin{multline}   \label{eq:inclchain}
 \bmond(\T,\LL(\h))
   \subsetneq \bmol(\T,\LL(\h)) \subsetneq \bmosd \\
       \subsetneq {\rm SBMO}(\T,\LL(\h))\subsetneq
{\rm WBMO}(\T,\LL(\h)).
\end{multline}

\lspace The reader is referred to  \cite{B1}, \cite{BP}, \cite{mei}, \cite{psm} for some
recent  results on dyadic BMO and Besov spaces connected to the
ones in this paper.

\lspace Mei's result implies in particular that
$\bmond(\T,\LL(\h))\nsubseteq \bmop$, and it is also easy to see that
the reverse inclusion does not hold (see for example the proof of $\bmol \nsubseteq \bmop$ at the beginning of Section
\ref{sec:2}).

To retrieve an estimate of the norm of the paraproduct in terms of the $\bmond$ norm,
we will consider the ``\emph{sweep}'', which is of independent interest, in Section \ref{sec:2},
and averages of martingale transforms in Section \ref{sec:4}.

Given $B\in \opf$, we define the sweep of $B$ as
\begin{equation} \label{sweep}
   S_B=\sum_{I \in \DD} B^*_I B_I\frac{\chi_I}{|I|}     .
\end{equation}
Our main result of Section \ref{sec:2}, Theorem \ref{mainteo}, states
that $$\|B\|^2_{\bmop}\approx \|S_B\|_{\bmol} + \|
B\|_{\sbmod}^2.$$ In particular, using the result
$\bmond(\T,\LL(\h))
   \subsetneq \bmol(\T,\LL(\h))$ (see \cite {new}), this shows that
if $B\in \sbmod$ and $S_B\in \bmond$, then $\pi_B$ is bounded.

Section 3 is devoted to the study of sweeps of functions in
different BMO-spaces. The classical John-Nirenberg theorem on
$\bmodd (\T)$ implies (and is essentially equivalent to) the fact
that there exists a constant $C >0$ such that
\begin{equation}\label{bmosweep}\|S_b\|_{\bmodd}\le C
\|b\|^2_{\bmodd}\end{equation} for any $b\in {\bmodd}$.

We will show that this formulation of John-Nirenberg does not hold
for $\|B\|_{\bmos}$. In fact, it is shown that if (\ref{bmosweep})
holds for some space contained in $\sbmod$ then this space is also
contained in $\bmop$.

In \cite{katz}, \cite{ntv} and \cite{nptv}, the correct rate of
growth of the constant in the Carleson embedding theorem in the
matrix case in terms of the dimension of Hilbert space $\h$ was
determined, namely $\log (\dim \h +1) $. Here, we want to show
that this  breakdown of the Carleson embedding theorem in the
operator case is intimately connected to a breakdown of the
John-Nirenberg Theorem, and that the dimensional growth for
constants in the John-Nirenberg Theorem is the same. This answers
a question left open in \cite{gptv2}.

In Section 4, we investigate ``average BMO conditions'' in the
following sense. We show  (see Theorem \ref{bmopara1}) that
$\|B\|_{\bmop}\le C(\int_\Sigma \|T_\sigma
B\|^2_{\bmond}d\sigma)^{1/2}.$ More precisely, $\|B\|^2_{\bmop}+
\| B^*\|_{\bmop}^2 \approx \int_\Sigma \|T_\sigma
B\|^2_{\bmol}d\sigma$.

Moreover, the norms $\|B\|_{\bmosd}$, $\|B\|_{\bmol}$
and $\|B\|_{\bmop}$ can be completely described in terms of average
boundedness
of
certain operators involving either
$\Lambda_B$ or commutators $[T_\sigma,B]$. The results of this section
complete
those proved
in \cite{gptv2}.

\section{Haar multipliers and paraproducts}   \label{sec:2}

We start by describing the action of a paraproduct $\pi_B$ as a Haar
multiplier.
\begin{prop}\label{equipara} Let $B\in \opf$. Then
$$\|\pi_B\|= \|(B^*_Ih_I)_{I\in \DD}\|_{mult}$$
$$=\|(P_{I^+}B+ P_{I^-}B)_{I\in \DD}\|_{mult}$$
$$=\|(\sum_{J\subsetneq I}B_J^* B_J\frac{\chi_J}{|J|})_{I\in
\DD}\|^{1/2}_{mult}.$$

In particular,
$$\|B_I\|\le \|\pi_B\||I|^{1/2},$$ $$\|P_{I^+}B(e)+
P_{I^-}B(e)\|_{L^2(\T,\h)}\le \|\pi_B\| |I|^{1/2}\|e\| $$ and  $$
\|(\sum_{J\subsetneq I}B_J^*
B_J\frac{\chi_J}{|J|})e\|_{{L^2(\T,\h)}} \le  \|\pi_B\|^2 |I|\|e\|
.$$
\end{prop}
\proof The first and second equalities follow directly from the
definitions and $\|\pi_B\|=\|\Delta_{B^*}\|$.

For the third equality, use $\|\pi_B\|^2=\|\pi_B^*\pi_B\|$,
\begin{eqnarray*}
\pi_B^*\pi_B(f)(t)&=&\sum_{I\in \DD} B_I^*B_I
(m_I(f))\frac{\chi_I(t)}{|I|}
=\sum_{I\in \DD} B_I^*B_I (\sum_{I\subsetneq J}f_J
m_I(h_J))\frac{\chi_I(t)}{|I|}\\
&=&\sum_{I\in \DD} B_I^*B_I (\sum_{I\subsetneq
J}f_J)h_J(t)\frac{\chi_I(t)}{|I|}
=\sum_{J\in \DD}  (\sum_{I\subsetneq
J}B_I^*B_I\frac{\chi_I(t)}{|I|})f_Jh_J(t).
\end{eqnarray*}

The estimates now follow  from (\ref{estimacion}). \qed

\lspace
 The following characterizations of  $\mathrm{SBMO}$ will be useful below.

\begin{prop} (\cite{gptv2}) \label{carbmoso} Let $B\in{\rm
SBMO^d}(\T,\LL(\h))$. Then

$$\|B\|^2_{{\rm SBMO^d}} =
\displaystyle\sup_{I\in \DD,\|e\|=1}\frac{1}{|I|}\|P_I(
B_e)\|^2_{L^2(\h)}  \approx\displaystyle\sup_{I\in
\DD}\frac{1}{|I|}\|\sum_{J\subseteq I} B_J^* B_J\|.
$$
\end{prop}

\lspace It follows at once from Propositions \ref{equipara}  and
\ref{carbmoso} that
$$\bmop (\T,\LL(\h))\subseteq \sbmod(\T,\LL(\h)).$$

It is easily seen that, if $B$ and $B^*$ belong to  $\bmop$, then $B\in
\bmol$.
However, we want to remark that the boundedness of $\pi_B$
alone does not
imply boundedness of $\Lambda_B$.

 To see this, choose some orthonormal basis $(e_i)_{i \in \N}$
of $\h$, and choose a sequence of
$\C^n$-valued function $(b_n)_{n \in \N}$ with finite Haar expansion
such that $\|b_n\|_{\bmoh} \ge C n^{1/2} \|b_n\|_{\bmow}$ (for a choice of
such a sequence,
see \cite{jpp1}).
Let $B_n(t)$ be the column matrix with respect to
the chosen orthonormal basis which has the
vector $b_n(t)$ as its first column. Then it is easy to see that
$$\|\pi_{B_n}\| = \| \pi_{b_n}\| \sim \|b_n\|_{\bmodd(\T,\HH)} \ge n^{1/2} C \| b_n
\|_{\wbmod(\T,\HH)}.$$
As pointed out to us \cite{pvpers},
it follows from the first Theorem in the appendix in \cite{pxu}
that
$\|\pi_{B_n^*}\| \le C \| b_n \|_{\wbmod(\T,\HH)}$ for some absolute constant $C$ and
all $n \in \N$.
Forming the direct sum
$$
  B = \bigoplus_{n=1}^\infty \frac{1}{\|\pi_{B_n^*}\|}B_n^*,
$$
we find that $\|\pi_B\|=1$, but $\Delta_B =(\pi_{B^*})^*$ is unbounded.

\lspace

\medskip

One of the main tools to investigate the connection between $\bmol$ and $\bmop$
is the
\emph{dyadic sweep}.
Given $B\in \opf$, we define
$$
S_B(t)=\sum_{I\in  \DD} B^*_IB_I
\frac{\chi_I(t)}{|I|}.
$$

\begin{lemm}\label{maincor} Let $B\in \opf$. Then
\begin{equation} \label{eq:sweepid}
     \pi_B^* \pi_B = \pi_{S_B} + \pi_{S_B}^* + D_B =\Lambda_{S_B}+D_B,
\end{equation}
where $D_B$ is defined by $D_B h_I \otimes x = h_I \frac{1}{|I|}
\sum_{J \subsetneq I} B_J^* B_J x$ for $x \in \h$, $I \in \DD$ and
$$ \|D_{B} \| \approx \|B\|^2_{\sbmod}. $$
\end{lemm}
\proof (\ref{eq:sweepid}) is verified on elementary tensors $h_I \otimes x$, $h_J \otimes
y$.
We find that
\begin{enumerate}
\item for $I \subsetneq J$,
$$\langle  \pi_B^* \pi_B h_I \otimes x, h_J \otimes y \rangle
 = \langle \pi_{S_B}^* h_I \otimes x, h_J \otimes y \rangle
$$
\item  for $I \supsetneq J$,
$$\langle  \pi_B^* \pi_B h_I \otimes x, h_J \otimes y \rangle
 = \langle \pi_{S_B} h_I \otimes x, h_J \otimes y \rangle
$$
\item for $I=J$,
$$\langle  \pi_B^* \pi_B h_I \otimes x, h_J \otimes y \rangle
 = \langle D_{B} (h_I \otimes x), h_J \otimes y \rangle.
$$
\end{enumerate}
Since $\supp \pi_{S_B} h_I \subseteq I$ and $\supp \Delta_{S_B} h_I
\subseteq
I$,
$\langle  \pi_B^* \pi_B h_I \otimes x, h_J \otimes y \rangle =0$ in all
other cases.

One sees easily that $D_{B}$ is block diagonal with respect
to the Hilbert space decomposition
$L^2(\T, \h) = \bigoplus_{I \in \DD} \h$ defined by the mapping
$f \mapsto (f_I)_{I \in \DD}$.
The operator $\pi_{S_B}$ is block-lower triangular with respect to
this decomposition (using the natural partial order on $\DD$), and
$\Delta_{S_B}$ is block-upper triangular. Thus we obtain the required
identity.
Note that
$$ \|D_{B} \| =\sup_{I\in \DD, \|e\|=1} \frac{1}{|I|}\|\sum_{J
\subsetneq I}
         B_J^* B_J e\| \approx \|B\|^2_{\sbmod}$$
by Proposition \ref{carbmoso}. \qed

Notice that $(S_B)^*=S_B$. Hence Lemma \ref{maincor} gives
\begin{satz}\label{mainteo}
$$\|S_B\|_{\bmol}+ \|B\|^2_{\sbmod} \approx  \|\pi_B\|^2.  $$
\end{satz}
         \proof
It suffices to use that $\|D_{B}\| \approx
\|B\|^2_{\sbmod}$ and that
$\|B\|_{\sbmod} \lesssim \|\pi_B\|$ (using Proposition \ref{equipara}).
\qed

\noindent
This provides, among other things, our first link between ${\bmond}$ and ${\bmop}$:
\begin{cor} \label{cor:paranorm}
$$
    \|\pi_B\|^2 \lesssim \|S_B\|_{\bmond} + \|B\|^2_{\bmosd}.
$$
\end{cor}
\proof Theorem \ref{mainteo} and (\ref{eq:inclchain}). \qed

\section{Sweeps of operator-valued functions}   \label{sec:3}

Let us mention that by John-Nirenberg's lemma,
 we actually have that $f\in \bmond$ if and only if
$$\sup_{I \in \DD}
(\frac{1}{|I|} \int_I \| B(t) - m_I B \|^p  dt)^{1/p}<\infty$$ for
some (or equivalently, for all) $0 < p < \infty$. Since $(B-
m_IB)\chi_I= P_IB$, we can also say that $f\in \bmond$ if and only
if
$$ \sup_{I \in \DD} \frac{1}{|I|^{1/p}}\|P_I( B) \|_{L^p(\LL(\h))}<
\infty.$$

  One way to express the John-Nirenberg
inequality on scalar-valued $\bmodd$ is to say that the mapping
\begin{equation}  \label{eq:sweepact}
    \bmodd \rightarrow \bmodd, \quad b \mapsto S_b,
\end{equation}
is bounded.
In the operator-valued setting, this John-Nirenberg property breaks
down. Our
main result is that any space of operator-valued functions which is
contained in
$\bmosd(\T, \LL(\h))$ and on which the mapping (\ref{eq:sweepact})
acts boundedly is already contained in $\bmop(\T,\LL(\h))$.

However, we find that (\ref{eq:sweepact}) acts boundedly between
different operator-valued BMO spaces. We also obtain the
sharp rate of growth of the norm of the mapping (\ref{eq:sweepact}) on
$\bmosd(\T, \LL(\h))$, $\bmop(\T, \LL(\h))$, $\bmol(\T, \LL(\h))$ and
$\bmond(\T, \LL(\h))$ in terms of the dimension of $\h$.

Before establishing this dimensional growth, we consider an extension of the sweep.
In the scalar case, one can extend the sweep
$\bmodd \rightarrow \bmodd$ to a sesquilinear map $\Delta: \bmodd \times \bmodd \rightarrow \bmodd$.
This map is motivated by the consideration of ``products of paraproducts'' $\pi_f^* \pi_g$,
which in turn is motivated by the long-standing investigation of products of Hankel operators
$\Gamma_f^* \Gamma_g$ in the literature (see \cite{psm} and the references therein).

\begin{defi} Let us denote by $\Delta: \opf\times \opf \to L^1(\T, \LL(\HH))$ the bilinear
map
given by
$$\Delta(B, F)=\sum_{I\in \DD} B^*_IF_I\frac{\chi_I}{|I|} .$$
In particular $ S_B =\Delta(B,B)$ and $\Delta(B, F)^*=\Delta(F, B)$.
\end{defi}

\begin{lemm}\label{comp} Let $B\in \opf$. Then
$$P_I\Delta(B, F)=P_I\Delta(B,P_I F)= P_I \sum_{J \subseteq I}
\frac{\chi_J}{|J|}
B_J^* F_J =  P_I \sum_{J \subsetneq I}
\frac{\chi_J}{|J|}
B_J^* F_J.$$
In particular, $P_I(S_B)=P_I(S_{P_IB})= P_I(S_{(P_{I^+} + P_{I^{-}}) B})$.
\end{lemm}
\proof $P_I\Delta({B^*},(F_Jh_J))= P_I(B^*_JF_J\frac{\chi_J}{|J|})=0$ if
$I\subseteq J$. Hence
$$ P_I\Delta(B, F)=P_I\Delta(B,P_I F) = P_I \Delta(B, (P_{I^+} +P_{ I^-}) F).$$
\qed

\lspace

A similar proof as in Lemma \ref{maincor}  shows that
\begin{lemm}\label{mainlema} Let $B, F\in \opf$. Then
\begin{equation*}
     \pi_B^* \pi_F = \pi_{\Delta(B,F)} + \pi_{\Delta(F,B)}^* + D_{B,F}
=\Lambda_{\Delta(B,F)}+D_{B,F},
\end{equation*}
where $D_{B,F}$ is defined by $D_{B,F} (h_I \otimes x) = h_I
\frac{1}{|I|} \sum_{J \subsetneq I} B_J^* F_J x$ for $x \in \h$, $I
\in \DD$.
 Moreover, $ \|D_{B,F} \| \le
\sup_{\|e\|=1}\|B_e\|_{BMO(\h)}\sup_{\|e\|=1}\|F_e\|_{BMO(\h)}.$
\end{lemm}

Let us now study the boundedness of the sesquilinear map $\Delta$ in the various
BMO norms. Again, the properties of the map $\Delta$ are more subtle in the operator-valued case
than in the scalar case.

\begin{satz} \label{maindelta} There exists a constant $C >0$ such that for
$B, F\in \opf$,

(i) $\|\Delta(B,F)\|_{\bmol}\le C \|B\|_{\bmop} \|F\|_{\bmop}, $

(ii) $\|\Delta(B,F)\|_{\wbmod}\le  C\|B\|_{\sbmod} \|F\|_{\sbmod}
$

(iii) $\|\Delta(B,F)\|_{\sbmod}\le  C\|\pi_B\| \|F\|_{\sbmod}. $

\end{satz}
\proof  (i) follows from Lemma \ref{mainlema}.

\lspace

\noindent
(ii) Using Lemma \ref{comp}, one obtains
$$\langle P_I\Delta(B, F)e, f\rangle= P_I \sum_{J\in \DD}\langle (P_IF)_Je,
(P_IB)_Jf\rangle\frac{\chi_J}{|J|}$$
for $e,f \in \HH$. Therefore,
\begin{eqnarray*}
\|\langle P_I\Delta(B, F)e, f\rangle\|_{L^1} &=&\|P_I \sum_{J\in
\DD}\langle (P_IF)_Je, (P_IB)_Jf\rangle\frac{\chi_J}{|J|}\|_{L^1}\\
&\le& 2 \|\sum_{J\in \DD}\langle (P_IF)_Je,
(P_IB)_Jf\rangle\frac{\chi_J}{|J|} \|_{L^1} \\
 &\le& 2\|(\sum_{J\in \DD}\|(P_IB)_Jf\|^2\frac{\chi_J}{|J|})^{1/2}\|_{L^2}
\|(\sum_{J\in \DD}\| (P_IF)_Je\|^2\frac{\chi_J}{|J|})^{1/2}\|_{L^2}\\
&\le&  2(\sum_{J\in \DD}\| (P_IB)_Jf\|^2)^{1/2}(\sum_{J\in \DD}\|
(P_IF)_Je\|^2)^{1/2}.
\end{eqnarray*}
Thus if $\|B\|_{\bmosd}= \|F\|_{\bmosd}=1$, then
$$\|\langle P_I\Delta(B, F)e, f\rangle\|_{L^1}\le
2\|P_IB_f\|_{L^2(\h)}\|P_IF_e\|_{L^2(\h)}\le 2|I|.$$ This, again
using John-Nirenberg's lemma, gives $\|\Delta(B,F)\|_{\bmow}\le C$.

\lspace

\noindent
(iii) From Lemma \ref{comp}, we obtain
$$\|P_I\Delta(B, F)e\|_{L^2(\h)}=\| \Delta_{B^*}(P_I F_e)\|_{L^2(\h)}\le
\|\pi_B\|\|P_I F_e\|_{L^2(\h)}.$$ \qed

\lspace

Here comes the main result of this section.
\begin{satz}   \label{thm:main5}
Let $\h$ be a separable, finite or infinite-dimensional Hilbert space.
Let $ \rho$ be a positive homogeneous functional  on the space $\opf$ of
$\LL(\h)$-valued functions
on $\T$
with finite formal Haar expansion such that there exists constants $c_1$,
$c_2$ with
\begin{enumerate}
\item $ \|B \|_{\bmosd} \le c_1 \rho(B)$ and
\item $  \rho(S_B) \le c_2 \rho(B)^2$ for all $B \in \opf$.
\end{enumerate}
Then there exists a constant $C$, depending only on $c_1$ and $c_2$, such
that
$\| B \|_{\bmop} \le C \rho(B)$ for all $B \in \opf$.
\end{satz}

\proof For $n \in \N$, let $E_n$ denote the subspace
$\{ f \in L^2( \T, \h): f_I =0 \text{ for } |I| < 2^{-n}\}$ of $L^2(\T,
\h)$.
Let
$c(n) = \sup\{ \|\pi_B\|_{E_n}: \rho(B) \le 1 \}$. An elementary
estimate
shows that $c(n)$ is well-defined and finite for each $n \in \N$.
For $\eps >0$, $n \in \N$, we can find $f \in E_n$, $\|f \|=1$, $B \in
\opf$,
$\rho(B) \le 1$  such that
\begin{multline*}
    c(n)^2 (1-\eps)^2 \le \| \pi_B f \|^2
            = \langle \pi_{S_B} f, f \rangle +
               \langle f, \pi_{S_B} f \rangle  +  \langle D_B f, f \rangle\\
        \le 2 c(n) \rho(S_B) + c_1 \|B \|_{\bmosd}
        \le 2 c_2 c(n) + c_1.
\end{multline*}
It follows that the sequence $(c(n))_{n \in \N}$ is bounded by $C = c_2 +
\sqrt{c_2^2 +c_1}$,
and therefore $\| \pi_B \| \le C \rho(B)$ for all $B \in \opf$.
\qed

One immediate consequence is the following answer to Question 5.1 in
\cite{gptv2}.
\begin{satz} There exists an absolute constant $C >0$ such that for
each $n \in \N$ and each measurable function $B: \T \rightarrow  \mat$,
\begin{equation}  \label{eq:sharpest}
    \| S_B \|_{\bmosd} \le C \log(n+1) \| B\|^2_{\bmosd},
\end{equation}
and this is sharp.
\end{satz}

\proof
From (iii) in Theorem \ref{maindelta} one obtains:
$$\|S_B\|_{\bmos}\le C\|B\|_{\bmop}\|B\|_{\bmosd}\le
C \log(n+1) \|B\|^2_{\bmosd},$$ since there exists an absolute
constant $C >0$ with
$$
   \|B\|_{\bmop} \le C \log(n+1) \|B\|_{\bmosd}
$$
by \cite{katz} and \cite{ntv}. On the other hand, denoting by $C_n$ the
smallest
constant such that
$$
 \| S_B \|_{\bmosd} \le C_n \| B\|^2_{\bmosd}
$$
for each  integrable function $B : \T \rightarrow
{\mathrm{Mat}}(\C, n \times n)$,
we obtain from Theorem \ref{thm:main5} that
$$
  \|B\|_{\bmop} \le (C_n + \sqrt{C_n^2 +1}) \|B \|_{\bmosd} \le 3 C_n  \|B
\|_{\bmosd}
$$
for each integrable $B$. It was shown in \cite{nptv} that
there exists an
absolute constant $c >0$ such that for each $ n \in \N$, there exists
$B^{(n)}: \T
\rightarrow \mathrm{Mat}(n \times n, \C)$ such that
$\|B^{(n)}\|_{\bmop} \ge c \log(n+1)  \|B^{(n)} \|_{\bmosd}$. Therefore
$C_n \ge \frac{c}{3} \log(n+1)$, and (\ref{eq:sharpest}) is sharp.
\qed

\lspace

Sharp rates of dimensional growth can also be determined for $S_B$ in $\bmond$, $\bmop$ and $\bmol$. Interestingly,
the rate of growth for $\bmosd$ and $\bmop$ is slower than the one for
$\bmol$ and $\bmond$.
\begin{satz} There exists an absolute constant $C >0$ such that for
each $n \in \N$ and each measurable function $B: \T \rightarrow  \mat$,
\begin{equation}  \label{eq:sharppara}
    \| S_B \|_{\bmop} \le C \log(n+1) \| B\|^2_{\bmop},
\end{equation}
\begin{equation}  \label{eq:sharpmult}
    \| S_B \|_{\bmol} \le C (\log(n+1))^2 \| B\|^2_{\bmol},
\end{equation}
\begin{equation}  \label{eq:sharpnorm}
    \| S_B \|_{\bmond} \le C (\log(n+1))^2 \| B\|^2_{\bmond},
\end{equation}
and this is sharp.

Corresponding estimates also hold for  the sesquilinear map
$\Delta$.
\end{satz}
\proof This is contained in \cite{new}. \qed

\lspace

\noindent

\noindent
Finally, the following corollary to Theorem \ref{thm:main5}   gives an estimate of $\|\cdot\|_{\bmop}$ in terms of
$\| \cdot\|_{\sbmod}$ with an ``imposed'' John-Nirenberg property. We need some notation:
Let $S^{(0)}_B = B$ and let $S^{(n)}_B = S_{S^{(n-1)}B}$ for $n \in \N$, $B \in \opf$.

\begin{cor} There exists a constant $C>0$ such that
$$
   \| B\|_{\bmop} \le C \sup_{n \ge 0} \|S^{(n)}_B\|^{1/2^n}_{\sbmod} \qquad (B\in \opf).
$$
\end{cor}

\proof Define $\rho(B) = \sup_{n \ge 0}
\|S^{(n)}_B\|_{\sbmod}^{1/2^n}$. One sees easily that this
expression is finite for $B \in \opf$. Now apply Theorem
\ref{thm:main5}. \qed

\section{Averages over martingale transforms and operator-valued  BMO} \label{sec:4}

Let $\Sigma = \{-1,1\}^\DD$, and let $d\sigma$ denote the natural
product probability measure on $\Sigma$, which assigns measure
$2^{-n}$ to cylinder sets of length $n$.

For $\sigma \in \{-1,1 \}^\DD$, define the \emph{dyadic martingale
transform}
\begin{equation}   \label{eq:martdef}
T_\sigma: L^2(\T, \h) \rightarrow L^2(\T, \h), \qquad f = \sum_{I
\in \DD} h_I f_I \mapsto \sum_{I \in \DD} h_I \sigma_I f_I,
\end{equation}
 Given a Banach space $X$ and $F\in L^1(\T, X)$, we write $\tilde
F$ for the function defined a.e. on $\Sigma\times \T$ by
$$ \tilde F(\sigma, t)= T_\sigma F(t)=\sum_{I}\sigma_I F_I h_I(t).$$

In case that $X$ is a Hilbert space, $\|T_\sigma F\|_{L^2(\T,
X)}=\|F\|_{L^2(\T, X)}$ for any $(\sigma_I)_{I \in \DD}$, and
therefore $\|\tilde F\|_{L^\infty(\Sigma, L^2(\T, X))}= \| F\|_{
L^2(\T, X)}$. More generally, we have for UMD spaces that
$\|T_\sigma F\|_{L^2(\T, X)}\approx\|F\|_{L^2(\T, X)}$. However,
$X=\LL(\h)$ is not a UMD space, unless $\h$ is finite dimensional.

Whilst $\|B\|_{\bmop}$ cannot be estimated in terms of
$\|B\|_{\bmol}$ \cite{mei}, we will prove an estimate of
 $\|B\|_{\bmop}$ in terms of an average
of $\|T_\sigma B\|_{\bmol}$ over $\Sigma$. Similarly, whilst the
result in \cite{mei} implies that $\|S_B\|_{\bmond}$ cannot be
estimated in terms of $\|B\|_{\bmond}$, we will prove an estimate
of
 $\|S_B\|_{\bmond}$ in terms of an average
of $\|T_\sigma B\|_{\bmond}$ over $\Sigma$. For this, the
following representation of the sweep will be useful:
\begin{equation} \label{eq:sweepav}
   S_B(t) = \int_{\Sigma} (T_\sigma B)^*(t) (T_\sigma B)(t) d\sigma.
\end{equation}

\begin{satz}\label{bmopara1} Let $B\in \opf$. Then
\begin{equation*}
\|S_B\|_{\bmond}\lesssim \int_\Sigma \|T_\sigma B
\|^2_{\bmond}d\sigma .
\end{equation*}

In particular $\|B\|_{\bmop}^2 \lesssim \int_\Sigma \|T_\sigma B
\|^2_{\bmond}d\sigma$.
\end{satz}
\proof The first inequality follows from the estimate
\begin{eqnarray*}
&&\|P_IS_B\|_{L^1(\T,\LL(\h))}
\\&=&\|P_IS_{P_IB}\|_{L^1(\T,\LL(\h))}
 \le2\left\|\int_\Sigma (T_\sigma P_IB^*)(T_\sigma P_IB)
d\sigma \right\|_{L^1(\T,\LL(\h))} \\
&\le&   2\int_\Sigma \|(P_IT_\sigma B)^*P_IT_\sigma
B\|_{L^1(\T,\LL(\h))} d\sigma
=2\int_\Sigma \|(P_IT_\sigma B)\|^2_{L^2(\T,\LL(\h))} d\sigma\\
&\le&2 |I|\int_\Sigma \|T_\sigma B \|^2_{\bmond} d\sigma.
\end{eqnarray*}
Using John-Nirenberg's lemma for $\bmond(\T,\LL(\h))$, one
concludes the result. The second inequality follows from the first, (\ref{eq:inclchain})
and Theorem \ref{mainteo}.\qed

\lspace

We are going to describe the different operator-valued BMO spaces
in terms of "average boundedness" of certain operators. First we
see that the  $\bmosd$-norm can be described by ``average
boundedness'' of $\Lambda_B$.

\begin{satz}\label{avebmoso} Let $B\in \opf$, and let $\Phi_B$ be the map
$$\Phi_B : L^2(\T, \h) \rightarrow L^2(\T \times \Sigma, \h),
   \quad f \mapsto \Lambda_B T_\sigma f.$$
Then
\begin{equation*} \| \Phi_B\|= \sup_{\|f\|_{L^2(\h)}=1}(\int_\Sigma
\|\Lambda_B(T_\sigma f)\|^2_{L^2(\T,\h)}d\sigma)^{1/2}=
\|B\|_{\sbmod}.
\end{equation*}
In particular, $\|B\|_{\bmos}= \| \Phi_B\|+ \| \Phi_{B^*}\|.$
\end{satz}
\proof  Since $\Lambda_B (T_\sigma f)=\sum_{I\in \DD}P_I(B) f_Ih_I
\sigma_I$, we have
\begin{eqnarray*}
&& \int_{\Sigma} \int_\T \| (\Phi_B f)(t, \sigma)\|^2 dt
d\sigma\\& =&
   \int_{\Sigma} \int_\T \| (\Lambda_B T_\sigma f)(t)\|^2 dt d\sigma
=  \sum_{I \in \DD}\|P_I(B) f_Ih_I\|^2_{L^2(\h)} \\
&=& \sum_{I \in \DD} \frac{1}{|I|}  \int_I \| (B(t) - m_I B)
(\frac{f_I}{\|f_I\|})\|^2\|f_I\|^2 dt
 \le  \sup_{\|e\|=1}\|B_e\|^2_{BMO(\h)} \sum_{J \in \DD} \|f_J\|^2.
\end{eqnarray*}
The reverse inequality follows by considering functions $f=h_I e$,
where $e \in \h$, $I \in \DD$. \qed

\lspace

We require a further technical lemma, which shows that the $L^2$
norm of $\tilde B f$ may be decomposed in a certain way.
\begin{lemm}\label{mainpara}
Let $B\in \opf$ and $f \in L^2(\T, \HH)$.  Write $Bf=\pi_B
f+\Delta_B f+\gamma_B f$. Then
\begin{multline*}
\| \tilde B f\|^2_{L^2(\Sigma\times \T, \h)} \\=
\int_\Sigma\|\pi_{ T_\sigma B}(f)\|^2_{L^2(\h)} d\sigma +
\int_\Sigma\|\Delta_{ T_\sigma B}(f)\|^2_{L^2(\h)}
d\sigma+\int_\Sigma\|\gamma_{ T_\sigma B}(f)\|^2_{L^2(\h)} d\sigma
\end{multline*}
and
\begin{equation}\label{multtilda}
\| \Lambda_{\tilde B} f\|^2_{L^2(\Sigma\times \T, \h)} =
\int_\Sigma\|\pi_{ T_\sigma B}(f)\|^2_{L^2(\h)} d\sigma +
\int_\Sigma\|\Delta_{ T_\sigma B}(f)\|^2_{L^2(\h)} d\sigma.
\end{equation}
\end{lemm}

\proof   Observe that $m_I(T_\sigma B) h_I=(\sum_{I\subsetneq
J}\sigma_J B_J h_J)h_I.$ Hence

$$\gamma_{T_\sigma B}(f)=\sum_{I\in \DD}m_I(T_\sigma
B)(f_I)h_I=\sum_{J\in \DD}\sigma_J B_J(\sum_{I\subsetneq J}
f_Ih_I)h_J. $$ This shows that
$$
\int_\T\int_\Sigma \langle \pi_{T_\sigma B} f, \gamma_{T_\sigma
B}g \rangle d\sigma dt= \sum_{I\in\DD}\int_I\langle  B_I m_I f ,
B_I(\sum_{J\subsetneq I} g_Jh_J) \rangle \frac{\chi_I}{|I|}dt=0;
$$

$$
\int_\T\int_\Sigma \langle \gamma_{T_\sigma B} f, \Delta_{T_\sigma
B}g \rangle d\sigma dt= \sum_{I\in\DD}\int_I\langle
B_I(\sum_{J\subsetneq I} f_Jh_J) , B_I g_I \rangle
\frac{h_I}{|I|}dt=0;
$$
$$
\int_\T\int_\Sigma \langle \pi_{T_\sigma B} f, \Delta_{T_\sigma
B}g \rangle d\sigma dt= \sum_{I\in\DD}\int_I\langle  B_I m_I f ,
B_I g_I  \rangle \frac{h_I}{|I|}dt=0.
$$
To finish the proof, simply expand $\| \tilde
B(f)\|^2_{L^2(\Sigma\times \T, \h)}$ and $\| \Lambda_{\tilde
B}(f)\|^2_{L^2(\Sigma\times \T, \h)}$. \qed

Here is our desired estimate of $\|B\|_{\bmop} + \|B^*\|_{\bmop}$
in terms of an average over $\|\tilde B\|_{\bmol}$.
\begin{cor} \label{cor:avchar}  Let $B \in \opf$. Then
$$\frac{1}{2}(\|\pi_B\|+ \|\Delta_B\|)  \le \| \tilde
B\|_{L^2(\Sigma, \bmol)}\le \|\pi_B\|+ \|\Delta_B\|.$$
\end{cor}
\proof To show the first estimate, it is sufficient to use
(\ref{multtilda}) in Lemma \ref{mainpara},  the identity
$\|\Delta_B\| = \|\pi_{B^*}\|$ and the invariance of the right
hand side under passing to the adjoint $B^*$.

For the reverse estimate, note that
$$
 \int_{\Sigma} \| \tilde B \|_{\bmol}^2 d \sigma
  \le \int_{\Sigma} ( \|\Delta_{T_\sigma B}\| + \|\pi_{T_\sigma B}\|)^2 d \sigma
   =( \|\Delta_{B}\| + \|\pi_{B}\|)^2.
$$
\qed

\lspace

\section{Acknowledgement}
We thank V.~Paulsen for a helpful discussion on operator space structures. We also thank
Tao Mei for his personal communication of a preliminary version of \cite{mei}.

\end{document}